\documentclass[a4paper,oneside,12pt]{article}
\usepackage[T2A]{fontenc}
\usepackage[utf8]{inputenc}
\usepackage{amsmath,amssymb}
\usepackage{expdlist}
\usepackage{textpos}
\usepackage{enumerate}
\usepackage{indentfirst}
\usepackage{url}
\usepackage{paralist}
\usepackage[linesnumbered,boxed]{algorithm2e}
\usepackage{graphicx}
\usepackage{nccmath}
\usepackage{array}
\usepackage[export]{adjustbox}
\usepackage{subcaption}
\usepackage{amsfonts}
\usepackage{resizegather}

\newcolumntype{L}{>{\centering\arraybackslash}m{3cm}}
\newcolumntype{P}[1]{>{\centering\arraybackslash}p{#1}}

\newcommand{\Cov}{\mathrm{Cov}}

\graphicspath{ {img/} }

\title{Relations between the random variable $w_x$ and the Dirichlet divisor problem}
\author{Dmitry Pyatin/ pyatin@petrsu.ru}
\date{sub. 03.07.2020,\\ rev. 25.09.2021}

\begin{document}

\maketitle

\noindent \textbf{Abstract} \\

We have developed a heuristic showing that in the Dirichlet divisor problem for \textit{almost all} $n \in \mathbb{N}^{+}$: $$ R(n) \leq O(\psi(n)n^{\frac{1}{4}}) $$ where
$$ R(n) = \Big\lvert \sum_{x=1}^{n}\Big\lfloor\frac{n}{x}\Big\rfloor - n\log{n} - (2\gamma-1)n \Big\rvert $$ and
$ \psi(n) $ - any positive function that increases unboundedly as $ n \to \infty $. The result is achieved under the hypothesis:
$$ \Big \{\frac{n}{x} \Big \} \sim w_x $$ where $ w_x $ is uniformly distributed over $ [0,1) $ random variable with a values set $ \{0, \frac {1} {x}, \ldots, \frac{x-1}{x} \} $ and the value accepting probability $ p = \frac{1}{x} $.\\

The paper concludes with a numerical argument in support of the hypothesis being true.
It is shown that the expectation: $$\mu_{1} \Big[\sum_{x=1}^{n}\Big(\frac{n}{x} - \frac{x-1}{2x}\Big)
\Big]= (2n+1)H_{\lfloor\sqrt{n}\rfloor} - \lfloor\sqrt{n}\rfloor^{2} - \lfloor\sqrt{n}\rfloor + C$$ has deviation from $D(n)$ is less than $R(n)$ in absolute value for all $n < 10^{5}$. \\

\noindent\textbf{Conventions}
\begin{itemize}
 \item[] $\{x\}$ -- fractional part of $x$;
 \item[] $[a,b]$ -- least common multiple of $a$ and $b$;
 \item[] $(a,b)$ -- greatest common divisor of $a$ and $b$;
 \item[] $\mu_k[f(x)]$ -- $k$-th central moment of $f(x)$;
 \item[] $\gamma$ -- the Euler-Mascheroni constant;
  \item[] \textit{almost all, almost everywhere} -- all elements of the set, except for a zero measure subset;
 \item[] $C$ -- some constant.
\end{itemize}

\section*{Introduction}

The Dirichlet divisor problem is to determine the lower bound for $\theta$ in the remainder estimate:
$$R(n) = \Big\lvert \sum_{x=1}^{n}\Big\lfloor\frac{n}{x}\Big\rfloor - n\log{n} - (2\gamma-1)n \Big\rvert = O(n^{\theta+\epsilon})$$ where $D(n)=\sum_{x=1}^{n}\Big\lfloor\frac{n}{x}\Big\rfloor$ -- divisor summatory function.

Using the hyperbola method Dirichlet showed \cite{dirihlet12}, that $\theta \leq \frac{1}{2}$.
G. Voronoi (1903) proved \cite{voronoi13}, that $\theta \leq \frac{1}{3}$. Further the result has improved repeatedly. H. Iwaniec and C. J. Mozzochi (1988) showed \cite{iwaniecmozzochi722}, that $\theta \leq 7/22$.
The best known result belongs to M. Huxley (2003), he established \cite{huxley131416}, that $\theta \leq \frac{131}{416}$.

In 1916, G. H. Hardy and independently E. Landau proved \cite{hardylandau14}, that $\theta \geq \frac{1}{4}$, therefore it has been established that: $$ \frac{1}{4} \leq \theta \leq \frac{131}{416} $$

It is believed that $\theta = \frac{1}{4}$.
In this paper, we show that \textit{almost everywhere} $\theta = \frac{1}{4}$ under some hypothesis.\\

\section*{Content}

It can be shown (see appendix) by using the result $\theta < \frac{1}{2}$ (\cite{voronoi13}-\cite{huxley131416}), the hyperbola method and equality:
\begin{equation}\label{eq1}
    \sum_{x=1}^{\sqrt{n}}\Big\{\frac{n}{x}\Big\} = C\lfloor\sqrt{n}\rfloor + g(\lfloor\sqrt{n}\rfloor)
\end{equation}
that for $n \to \infty$ the constant $C=\frac{1}{2}$ and hold:
\begin{equation}\label{eq2}
    R(n) =  \Big\lvert \sum_{x=1}^{\sqrt{n}}\Big(\Big\{\frac{n}{x}\Big\} - \frac{1}{2}\Big)  \Big\rvert 
\end{equation}

From the work of J. Kubilius it is known \cite{kubilius}, that as $n \to \infty$:
\begin{equation}\label{eq3}
    \nu_n\{ \lvert R(n) - \mu_1[R(n)] \rvert \leq \psi(n)\sqrt{\mu_2[R(n)]}\} \to 1
\end{equation}, where $\nu_n\{\ldots\}$ -- is the frequency of events with a condition $\{\ldots\}$.
$\nu_n= \frac{1}{n}N_n\{\ldots\}$, where $N_n$-- number of events with a condition $\{\ldots\}$, and
$\psi(n)$ -- any positive function that increases unboundedly as $n\to\infty$.\\

Thus, for \textit{almost all} $n$ the following inequation is hold:
\begin{equation}\label{eq4}
    \lvert R(n) - \mu_1[R(n)] \rvert  \leq \psi(n)\sqrt{\mu_2[R(n)]}
\end{equation}

To find $\mu_1[R(n)]$ and $\mu_2[R(n)]$, we hypothesize:\\

\textbf{Hypothesis 1}
For an arbitrary positive integer $n$ and a fixed positive integer $x$:
\begin{equation}\label{eq5}
    \Big\{\frac{n}{x}\Big\} \sim w_x
\end{equation}, where $w_x$ is a random variable uniformly distributed on $[0,1)$, taking values from $\{0,\frac{1}{x},\ldots,\frac{x-1}{x} \}$ with probability $p=\frac{1}{x}$.\\

Accepting this hypothesis, we can find $\mu_1[R(n)]$:
\begin{equation}\label{eq6}
\begin{aligned}
    &\mu_1[R(n)] = \mu_1\Big[\sum_{x=1}^{\sqrt{n}} \Big( \Big\{ \frac{n}{x} \Big\} - \frac{1}{2} \Big)\Big] = \sum_{x=1}^{\sqrt{n}} \mu_1\Big[\Big\{ \frac{n}{x} \Big\} - \frac{1}{2}\Big] = \\
    &= \sum_{x=1}^{\sqrt{n}}\Big( \frac{1}{x} \sum_{k=0}^{x-1}\frac{k}{x} - \frac{1}{2}\Big) = -\sum_{x=1}^{\sqrt{n}}\frac{1}{2x} = -\frac{1}{2} H_{\lfloor\sqrt{n}\rfloor}
\end{aligned}
\end{equation}

and $\mu_2[R(n)]$:
\begin{equation}\label{eq7}
    \mu_2[R(n)] = \mu_2\Big[\sum_{x=1}^{\sqrt{n}}\Big(w_x-\frac{1}{2}\Big)\Big] = \sum_{x=1}^{\sqrt{n}}\sum_{y=1}^{\sqrt{n}}\Cov(w_x-\frac{1}{2},w_y-\frac{1}{2})
\end{equation}
where:
\begin{equation}\label{eq8}
\begin{aligned}
    &\Cov(w_a-\frac{1}{2},w_b-\frac{1}{2}) = \\ &= \frac{1}{[a,b]}\sum_{i=1}^{\frac{a}{(a,b)}}\sum_{j=1}^{\frac{b}{(a,b)}}\sum_{k=1}^{(a,b)}\Big( \frac{a - ((i-1)(a,b) + k)}{a} - \frac{a-1}{2a}\Big)\Big(\frac{b-((j-1)(a,b)+k)}{b} - \frac{b-1}{2b} \Big) = \\
    &= \frac{1}{[a,b]} \frac{(a,b)^2 - 1}{12(a,b)} = \frac{(a,b)}{12[a,b]} - \frac{1}{12ab}
\end{aligned}
\end{equation}

The formula (\ref{eq8}) comes from the block structure of the covariance matrix $A(d_1,d_2)$:
\begin{equation}\label{eq9}
A\left(\frac{a}{(a,b)},\frac{b}{(a,b)}\right) = \begin{bmatrix} 
    G_{11} & G_{12} & \dots \\
    \vdots & \ddots & \\
    G_{\frac{a}{(a,b)},1} &        & G_{\frac{a}{(a,b)},\frac{b}{(a,b)}} 
    \end{bmatrix}
\end{equation}
where $G(d_1,d_2)$ diagonal matrix:
\begin{equation}\label{eq10}
G((a,b),(a,b)) = \begin{bmatrix} 
    \frac{1}{[a,b]} & 0 & \dots \\
    \vdots & \ddots & \\
    0 &        & \frac{1}{[a,b]} 
    \end{bmatrix}
\end{equation}
using the general formula for finding the covariance of two discrete random variables:
\begin{equation}
    \Cov(w_a,w_b) = \sum_{i=1}^{a}\sum_{j=1}^{b}p_{i j}(w_i - \mu_1(w_a))(w_j - \mu_1(w_b))
\end{equation}
in which some terms are equal to zero due to $p_{i j}=0$, and the number of nonzero terms is equal to $[a,b]$.\\

After all we have a second central moment:
\begin{equation}\label{eq11}
    \mu_2[R(n)] = \sum_{x=1}^{\sqrt{n}}\sum_{y=1}^{\sqrt{n}}\Big(\frac{(a,b)}{12[a,b]} + \frac{1}{12ab}\Big) = \sum_{x=1}^{\sqrt{n}}\sum_{y=1}^{\sqrt{n}}\frac{(a,b)}{12[a,b]} + O((\log{\sqrt{n}})^2)
\end{equation}

L. Toth (et al.) gives \cite{toth} an explicit formula for the sum with GCD and LCM:
\begin{equation}\label{eq12}
    \sum_{a=1}^{n}\sum_{b=1}^{n} \frac{(a,b)}{[a,b]} = 3n + O((\log{n})^2)
\end{equation}
so:
\begin{equation}\label{eq13}
    \mu_2[R(n)] = \frac{1}{4}\lfloor\sqrt{n}\rfloor + O((\log{n})^2)
\end{equation}

whence by substitution in (\ref{eq3}) and adding $\frac{1}{2}H_{\lfloor\sqrt{n}\rfloor}$ we get that for $n\to \infty$:
\begin{equation}\label{eq14}
    \nu_n\{ R(n) > O(\psi(n)n^{\frac{1}{4}})\} \to 0
\end{equation} where $\psi(n)$ -- any positive function, increasing unboundedly as $n\to\infty$, and as a consequence as $n\to \infty$ \textit{almost everywhere} performed:
\begin{equation}
    R(n) \leq O(\psi(n)n^{\frac{1}{4}})
\end{equation}

\section*{Conclusion}

Let's try to present an argument in support of the truth of hypothesis 1.
Knowing that:
\begin{equation}\label{eq15}
    D(n) = \sum_{x=1}^{n}\Big\lfloor\frac{n}{x}\Big\rfloor = \sum_{x=1}^{n}\frac{n}{x} - \sum_{x=1}^{n}\Big\{\frac{n}{x}\Big\}
\end{equation}
Let's define a random variable:
\begin{equation}\label{eq16}
    W(n) = \sum_{x=1}^{n}\frac{n}{x} - \sum_{x=1}^{n}w_x
\end{equation}

Find $\mu_1[W(n)]$, using $\mu_1[w_x] = \frac{x-1}{2x}$ and the Dirichlet hyperbola method:
\begin{equation}\label{eq17}
    \mu_1[W(n)] = (2n+1)H_{\lfloor\sqrt{n}\rfloor} - \lfloor\sqrt{n}\rfloor^{2} - \lfloor\sqrt{n}\rfloor + C
\end{equation}

Numerical calculations show that $\mu_1[W(n)]$ closer to $D(n)$, that $n\log{n} + (2\gamma-1)n$. Introduce the error functions:
\begin{equation}
\begin{aligned}
    \Delta_R &= \sum_{n=1}^{N} \lvert R(n) \rvert, \;
    \Delta_W = \sum_{n=1}^{N} \Big(\lvert D(n) - \mu_1[W(n)]\rvert\Big) \\
    d_R &= \sum_{n=1}^{N} R(n), \;
    d_W = \sum_{n=1}^{N} \Big(D(n) - \mu_1[W(n)]\Big)
\end{aligned}
\end{equation}

\begin{figure}[h!]
  \center
  \includegraphics[scale=0.55]{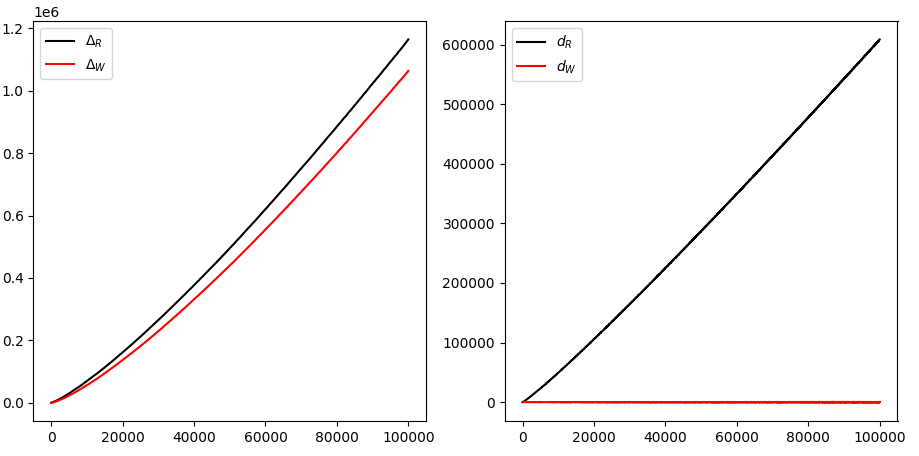}
  \caption{\label{img:pic1} Comparison of $\Delta_R$, $\Delta_W$, $d_R$ and $d_W$.}
\end{figure}

\begin{figure}[h!]
  \center
  \includegraphics[scale=0.6]{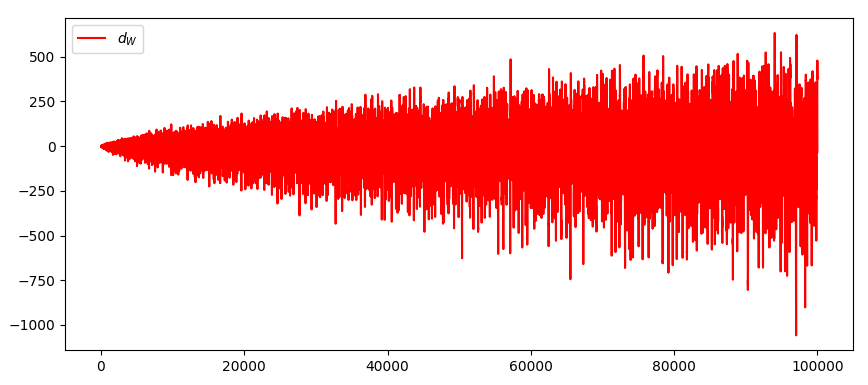}
  \caption{\label{img:pic2} $d_W$ up to $N=10^5$.}
\end{figure}

Figure \ref{img:pic1} shows that $\Delta_W < \Delta_R$ and $d_W \ll d_R$. Figure \ref{img:pic2} shows $d_W$ using the constant $C=\frac{1}{12}$ in the formula (\ref{eq17}). It can be replaced $\frac{1}{2}$ with $\frac{x-1}{2x}$ in formula (\ref{eq2} ) for getting $\mu_1 = 0$.\\

As a result the proposed heuristic is in good agreement with the numerical data for the parameter up to $10^5$, however, the proof of the estimates for $\theta$ in the Dirihlet divisor problem must be carried out using other methods.\\

\newpage
\section*{Appendix}
Let us prove that:
 \begin{equation}
    \sum_{x=1}^{\sqrt{n}} \Big\{ \frac{n}{x}\Big\} = \frac{1}{2}\lfloor\sqrt{n}\rfloor + R(n)
 \end{equation}
Using the equality:
 \begin{equation}
 \sum_{x=1}^{\sqrt{n}}\Big\{ \frac{n}{x}\Big\} = C \lfloor \sqrt{n} \rfloor + g(\lfloor\sqrt{n}\rfloor)
 \end{equation}
 using the Dirichlet hyperbola method, we obtain:
 \begin{equation}
   \begin{aligned}
     D(n) &= \sum_{x=1}^{n}\Big\lfloor \frac{n}{x}\Big\rfloor = 2\sum_{x=1}^{\sqrt{n}}\Big\lfloor \frac{n}{x}\Big\rfloor - \lfloor\sqrt{n}\rfloor^2 = 2\sum_{x=1}^{\sqrt{n}}\Big( \frac{n}{x} -\Big\{ \frac{n}{x}\Big\}\Big) - \lfloor\sqrt{n}\rfloor^2 = \\
     &= 2n\Big(\log{\lfloor\sqrt{n}\rfloor}) + \gamma + \frac{1}{2\lfloor\sqrt{n}\rfloor} + O(\frac{1}{n}\Big) - 2C\lfloor\sqrt{n}\rfloor + g(\lfloor\sqrt{n}\rfloor) - \lfloor\sqrt{n}\rfloor^2 = \\
     &= 2n\Big( \log{\sqrt{n} - \{\sqrt{n}\}} + \gamma + \frac{1}{2\lfloor\sqrt{n}\rfloor} + O(\frac{1}{n})\Big) - 2C\lfloor\sqrt{n}\rfloor + g(\lfloor\sqrt{n}\rfloor) - (\sqrt{n} - \{\sqrt{n}\})^2
    \end{aligned}
 \end{equation}
because:
 \begin{equation}
    \begin{aligned}
     &\log{\sqrt{n} - \{\sqrt{n}\}} = \log{\sqrt{n}}-\frac{\{\sqrt{n}\}}{\sqrt{n}} + O(\frac{\{\sqrt{n}\}}{n}), \text{and} \\
     & (\sqrt{n} - \{\sqrt{n}\})^2 = n - 2\{\sqrt{n}\}\sqrt{n} + \{\sqrt{n}\}^2
     \end{aligned}
 \end{equation}
 then, taking into account $\{\sqrt{n}\} < 1$:
 \begin{equation}
    \begin{aligned}
     D(n) &= 2n\log{\sqrt{n}} - \frac{2n\{\sqrt{n}\}}{\sqrt{n}} + O(\{\sqrt{n}\}) + 2\gamma n + \frac{2n}{2\lfloor\sqrt{n}\rfloor} + O(1) - \\
     &- 2C\lfloor\sqrt{n}\rfloor + g(\lfloor\sqrt{n}\rfloor) - n + 2\{\sqrt{n}\}\sqrt{n} - \{\sqrt{n}\}^2 \leq \\
     &\leq n\log{n} - 2\sqrt{n}\{\sqrt{n}\} + 2\{\sqrt{n}\}\sqrt{n} + 2\gamma n + \frac{n}{\lfloor\sqrt{n}\rfloor} - 2C\lfloor\sqrt{n}\rfloor - n + O(1)
    \end{aligned}
 \end{equation}
 because $\frac{n}{\lfloor\sqrt{n}\rfloor} < \lfloor\sqrt{n}\rfloor + 1$, so:
 \begin{equation}
     \begin{aligned}
        D(n) < n\log{n} + (2\gamma-1)n + (1-2C)\lfloor\sqrt{n}\rfloor + g(\lfloor\sqrt{n}\rfloor) + O(1)
     \end{aligned}
 \end{equation}
whence follows:
 \begin{equation}
     R(n) = D(n)-(n\log{n} + (2\gamma-1)n) = (1-2C)\lfloor\sqrt{n}\rfloor + g(\lfloor\sqrt{n}\rfloor)
 \end{equation}
т. к. $R(n) \ll \sqrt{n}$, under \cite{voronoi13} - \cite{huxley131416}, we conclude that: $$ C = \frac{1}{2} $$
As a result:
\begin{equation}
    R(n) = \sum_{x=1}^{\sqrt{n}} \Big(\Big\{ \frac{n}{x}\Big\} - \frac{1}{2}\Big) = -\frac{1}{\pi} \sum_{x=1}^{\sqrt{n}} \sum_{k=1}^{\infty}\frac{\sin{(2\pi k \frac{n}{x})}}{k}
\end{equation}
\end{document}